\documentstyle{amsppt}
\TagsOnRight \headline={\tenrm\hss\folio\hss}
\magnification=\magstep1\nopagenumbers \parindent=24pt\vsize 22
true cm \hsize 16 true cm

\topmatter
\title On the finitely separability of subgroups of generalized free products\endtitle
\author David Moldavanskii, Anastasiya Uskova  \endauthor

\abstract It is proved that all finitely generated subgroups of
generalized free product of two groups are finitely separable
provided that free factors have this property and amalgamated
subgroups are normal in corresponding factors and  satisfy the
maximum condition for subgroups.\endabstract
\endtopmatter
\document

\centerline{\bf 1. Introduction and results}
\medskip

Recall (see [1]), that subgroup $H$ of group $G$ is said to be
finitely separable if for any element $g\in G\setminus H$ there
exists a homomorphism $\varphi$ of group $G$ onto some finite
group such that the image $g\varphi$ of element $g$ does not
belong to the image $H\varphi$ of subgroup $H$. It is easy to see
that subgroup $H$ of group $G$ is finitely separable if and only
if $H$ coincides with the intersection of all subgroups of finite
index in $G$ containing $H$.

Hence it is obvious that any group $G$ is residual finite if and
only if the identity subgroup of $G$ is finitely separable.
Moreover, a normal subgroup $H$ of group $G$ is finitely separable
if and only if the quotient group $G/H$ is residual finite. Since
there exist 2-generated groups that are not residual finite (see
e.~g. [2]), this remark implies  that any non-cyclic free group
contains a subgroup which is not finitely separable. This, in
turn, implies that the problem of finding conditions under which
the free constructions of groups inherit the property of finitely
separability of all subgroups actually comes down to the question
of lack in the resulting groups of non-cyclic free subgroups. The
corresponding answers are well known.

On the other hand, from the theorem of M.~Hall [3] it follows that
every finitely generated subgroup of any free group is finitely
separable. Thus, an arbitrary free group is a subgroup separable
in the sense of the following definition currently used in a
number of publications: a group $G$ is called subgroup separable
if all finitely generated subgroups of $G$ are finitely separable.

Investigation of the conditions under which the groups obtained by
a free constructions are subgroup separable appeared to be more
interesting task. So, N.~Romanovskii [4] have proved that the
(ordinary) free product of any family of subgroup separable groups
is also a subgroup separable group.

It turned out, however, that for construction of generalized free
product of groups a similar statement is not true in general.
Apparently, the first (and the simplest) example of the
generalized free product of two subgroup separable groups which is
not a subgroup separable group was given in the work of R.~Allenby
and R.~Gregorac [5]. Authors of this work observed that 1) the
direct product $F_{2}\times F_{2}$ of two free groups of rank 2
contains a finitely generated subgroup that is not finitely
separable, 2) the direct product  $F_{2}\times F_{1}$ of free
group of rank 2 and of infinite  cyclic group is a subgroup
separable group and 3) the group $F_{2}\times F_{2}$ is
generalized  free product of two copies of the group $F_{2}\times
F_{1}$ with amalgamation $F_{2}$. An example of the generalized
free product with cyclic amalgamation of two subgroup separable
groups, which is not a subgroup separable group, was built by
E.~Rips [6]. Later the more simple similar example of generalized
free product with cyclic amalgamation of two finitely generated
nilpotent groups was given in [7].

Availability of these and other similar examples justifies
the interest to finding conditions  under which the
generalized free product $G=\bigl(A\ast B;\ H=K, \varphi\bigr)$ of
groups $A$ and $B$ with subgroups $H\leq A$ and $K\leq B$
amalgamated in accordance with isomorphism $\varphi:H\to K$
 inherit from factors the property of being a subgroup separable group.
For example, it is well known that if groups $A$ and $B$ are
finite then the group $G$ is subgroup separable (because in the
case when $A$ and $B$ are finite the group $G$ is free-by-finite
(see., e.~g., [8]) and as any finite extension of subgroup
separable group is a subgroup separable group too).  Let us recall
a few more similar results (the proofs of which are not, of
course, just as simple).

{\it The group $G=\bigl(A\ast B;\ H=K, \varphi\bigr)$ is subgroup
separable if one of the following  assumptions is
fulfilled:\roster

\item groups $A$ and $B$ are subgroup
separable and subgroups $H$ and $K$ are finite $[5]$;

\item groups $A$ and $B$ are polycyclic-by-finite and there exists
a finite index subgroup $U$ of group $H$ such that  $U$ and
$U\varphi$ are normal subgroups of groups $A$ and $B$ respectively
$[5]$;

\item groups $A$ and $B$ are free-by-finite and subgroups $H$ and
$K$ are cyclic $[9]$.\endroster}

The purpose of this article is to obtain a further result of a
similar nature.

 \proclaim{\indent Theorem} Let
$G=\bigl(A\ast B;\ H=K, \varphi\bigr)$ be the generalized free
product of groups $A$ and $B$ with subgroups $H\leq A$ and $K\leq
B$ amalgamated in accordance with isomorphism $\varphi:H\to K$.
Let $H$ be a normal subgroup of group $A$, $K$ be a normal
subgroup of group $B$ and groups $H$ and $K$ satisfy the maximum
condition for subgroups. If groups $A$ and $B$ are subgroup
separable then $G$ is a subgroup separable group.\endproclaim

Obviously, in the case when the groups $A$ and $B$  are
polycyclic-by-finite, the statement of the Theorem is contained in
the above result (2) of the work [5]. Let us state two another
corollaries of this Theorem assertions of which are, apparently,
the new. The first one is interesting to compare with the above
examples of [6] and [7] and result of [9].

\proclaim{\indent  Corollary 1} If $A$ and $B$ are subgroup
separated groups and   $H$ and $K$ are cyclic normal subgroups of
groups $A$ and $B$ respectively then the group  $G=\bigl(A\ast B;\
H=K, \varphi\bigr)$ is  subgroup separable.\endproclaim

In the second consequence  the question of subgroup separability
of Baumslag -- Solitar groups is considered.

Recall that Baumslag -- Solitar group is an one-relator group with
presentation
$$
G(m,n)=\langle a,\,b;\ a^{-1}b^{m}a=b^{n}\rangle,
$$
where $m$ and $n$ are non-zero integers. Note that since groups
$G(m,n)$, $G(n,m)$ and $G(-m,-n)$ are mutually isomorphic we can
assume without loss of generality that integers $m$ and $n$ in the
presentation of group $G(m,n)$  satisfy the condition
$|n|\geqslant m>0$. It is known (see [10]) that the group $G(m,n)$
is residual finite if and only if (under the condition
$|n|\geqslant m>0$) either $m=1$ or $|n|=m$.

It is also well known (and easily to see) that if  $|n|>1$ then in
group $G(1,n)$ the cyclic subgroup generated by element  $b$ is
not finitely separable. On the other hand if $|n|=m$ then the
group
$$
G(m,n)=\langle a,\,b;\ a^{-1}b^{m}a=b^{\pm m}\rangle=\langle
a,\,b,\,c;\ a^{-1}ca=c^{\pm 1},\ c=b^{m}\rangle
$$
is the generalized free product with normal amalgamated subgroups
of polycyclic group $G(1,\pm 1 )$ and infinite cyclic group. Thus,
we have

\proclaim{\indent Corollary 2} If $|n|=m$ then the group  $
G(m,n)=\langle a,\,b;\ a^{-1}b^{m}a=b^{n}\rangle$ is subgroup
separable.\endproclaim

It should be noted that the statement of Corollary 2 is new only
in the case when $n=-m$. Indeed, when $n=m$ then the center of
group $G(m,n)$ is non-trivial, and, as was proved in [11],  any
one-relator group with non-trivial center is subgroup separable.
\bigskip

\centerline{\bf 2. The proof of Theorem}
\medskip

We begin with the following simple remark.

\proclaim{\indent Lemma} Let $N$ be a normal subgroup of group
$L$. If $N$ is a finitely generated group and group $L$ is
subgroup separable then the quotient group $L/N$ is a subgroup
separable group.\endproclaim

Indeed, an  arbitrary subgroup of quotient group $L/N$ is
represented in the form $M/N$ for some subgroup $M$ of group $L$
containing $N$. In addition, as subgroup $M$ is generated by union
of system of representatives of cosets generating subgroup $M/N$
and of system of generators of subgroup $N$, if  subgroup $M/N$ is
finitely generated then subgroup $M$ is finitely generated too.

Now, let subgroup $M/N$ of quotient group $L/N$ be finitely
generated and let $aN$ be element of $L/N$ that does not belong to
subgroup $M/N$. Then in group $L$ element $a$ does not contained
in subgroup $M$ and since this subgroup is finitely generated and
therefore is finitely separable in group $L$ there exists a finite
index subgroup $R\leq L$ containing subgroup $M$ and not
containing element $a$. Then subgroup $R/N$ of quotient group
$L/N$ is of finite index, contains subgroup $M/N$ and not contains
element $aN$. Thus, subgroup $M/N$ coincides with the intersection
of all containing it subgroups of finite index of group $L/N$  and
therefore is finitely separable.\medskip

Suppose now that the group $G=\bigl(A\ast B;\ H=K, \varphi\bigr)$
satisfies all the assumptions in Theorem and therefore, in
particular, subgroup $H=K$ is normal in  $G$. Let $U$ be a
finitely generated subgroup of $G$ and let $a$ be an element of
$G$ not belonging to subgroup $U$. To prove the existence of the
homomorphism $\theta$ of group $G$ onto finite group such that
$a\theta\notin U\theta$ we consider separately the two cases.

{\it Case} 1. Element $a$ does not belong to subgroup $UH$. It is
easy to see that the quotient group $G/H$ is isomorphic to the
free product of groups $A/H$ and $B/K$. Since by Lemma these
groups are subgroup separable the result of [4] implies that $G/H$
is a subgroup separable group. As $UH/H\simeq U/(U\cap H)$,
subgroup $UH/H$ of group $G/H$ is finitely generated and since
$aH\notin UH/H$ there exists a homomorphism $\rho$ of group $G/H$
onto some finite group such that $(aH)\rho\notin (UH/H)\rho$. Then
the product $\theta$ of  the natural mapping of $G$ onto the
quotient group $ G/H$ and of homomorphism $\rho$ is the desired
homomorphism.\smallskip

{\it Case} 2. $a\in UH$ and therefore $a=uh$ for suitable elements
$u\in U$ and $h\in H$. As $a\notin U$ element $h$ does not belong
to subgroup $V=U\cap H$. Since group $H$  satisfies the maximum
condition, subgroup $V$ is finitely generated and therefore is a
finitely separable subgroup of group $A$. Consequently, there
exists a finite index normal subgroup $R$ of group $A$ such that
 $h\notin VR$. Then $T=H\cap R$ is a subgroup of finite index of group
$H$ and $h\notin VT$. Since group $H$ is finitely generated,
subgroup $T$ without loss of generality may supposed to be
characteristic in $H$ and therefore normal in $A$. Then
$S=T\varphi$ is characteristic in $K$ and normal in $B$.

We set,further, $\overline A=A/T$, $\overline H=H/T$, $\overline
B=B/S$ and $\overline K=K/S$. It is clear that the mapping
$\overline\varphi:\overline H\to\overline K$, defined by the rule
$(xT)\overline\varphi=(x\varphi)S$ for any $x\in H$, is the
isomorphism of group $\overline H$ onto group $\overline K$.
Therefore, we can construct the generalized free product
$\overline G=\bigl(\overline A*\overline B;\ \overline H=\overline
K, \overline \varphi\bigr)$ of groups $\overline A$ and $\overline
B$ with subgroups $\overline H$ and $\overline K$ amalgamated in
accordance with isomorphism $\overline \varphi$.

Since in group $\overline G$ for any $x\in H$  the equality
$xT=(x\varphi)S$ holds, the natural mappings of groups $A$ and $B$
onto subgroups $\overline A$ and  $\overline B$ of $\overline G$
are consistent with isomorphism $\varphi$. Therefore, there exists
homomorphism $\rho$ of group $G$ onto group $\overline G$
extending these mappings. It is easy to see that the kernel of
$\rho$ coincides with the subgroup $T$. Consequently, the
inclusion $a\rho\in U\rho$ is equivalent to the inclusion $a\in
UT$ and therefore is equivalent to the equality $a=u_{1}t$ for
suitable $u_{1}\in U$ and $t\in T$. Then $uh=u_{1}t$ whence we
obtain $u^{-1}u_{1}=ht^{-1}\in U\cap H=V$ and $h=(u^{-1}u_{1})t\in
VT$ which is impossible.

Thus, in group $\overline G$ element $a\rho$ does not belong to
finitely generated subgroup $U\rho$. By Lemma groups $\overline A$
and $\overline B$ are subgroup separable and since subgroups
$\overline H$ and $\overline K$ are finite the group $\overline G$
is subgroup separable by above result (2) of [5].  The proof of
the existence of the required homomorphism $\theta $ can now be
completed as in case 1. \newpage

\centerline{\bf References}
\medskip

\item{1.}
{\it Mal'cev A.~I.} On homomorphism onto finite groups
// Uchen. Zap. Ivanov. Ped. Inst. 1958. Vol.~18. P.~49--60
(Russian).

\item{2.}
{\it Baumslag G., Solitar D.} Some two-generator one-relator
non-Hop\-fian groups // Bull. Amer. Math. Soc. 1962. Vol.~68.
P.~199--201.

\item{3.}
{\it  Hall M.} Coset representations in free groups // Trans.
Amer. Math. Soc. 1949. Vol.~67. P.~421--432.

\item{4.}
{\it Romanovskii N.~S.} On the residual finiteness of free
products with respect to subgroups // Izv. Akad. Nauk SSSR Ser.
Mat. 1969. Vol.~33. P.~1324--1329 (Russian).

\item{5.}
{\it Allenby R., Gregorac R.} On locally extended residually
finite groups // Lecture Notes Math. 1973. Vol.~319. P.~9--17.

\item{6.}
{\it Rips E.} An example of a non-LERF group which is a free
product of LERF groups with an amalgamated cyclic subgroup //
Israel J. of Math., Vol.~70, ü~1, 1990. P.~104--110.

\item{7.}
{\it Allenby R., Doniz D.} A free product of finitely generated
nilpotent groups amalgamating a cycle that is not subgroup
separable // Proc. Amer. Math. Soc. Vol.~124, ü~4,1996.
P.~1003-1005.

\item{8.}
{\it Baumslag G.} On the residual finiteness of generalized free
products of nilpotent groups // Trans. Amer. Math. Soc. 1963.
Vol.~106. ü~2. P.~193--209.

\item{9.}
{\it Allenby R., Tang C.} Subgroup separability of generalized
free products of free-by-finite groups // Canad. Math. Bull. 1993.
Vol.~36(4). P.~385--389.

\item{10.}
{\it Meskin S.} Nonresidually finite one-relator groups // Trans.
Amer. Math. Soc. 1972. Vol.~164. P.~105--114.

\item{11.}
{\it Moldavanskii D.~I., Timofeeva L.~V.} Finitely generated
subgroups of one-relator group with non-trivial center are
finitely separable // Izvestiya Vuzov. Mathematics. 1987.
Issue~12. P.~58--59 (Russian).\bigskip

Ivanovo State University \smallskip

{\it E-mail address}: moldav\@mail.ru

\end